\documentclass[12pt]{amsart}

\usepackage[english]{babel}
\usepackage[utf8]{inputenc}
\usepackage{amsmath}
\usepackage{amssymb}
\usepackage{amsfonts}
\usepackage{amsthm}
\usepackage{mathrsfs}
\usepackage[all]{xy}
\usepackage[pdftex]{graphicx}
\usepackage{color}
\usepackage{cite}
\usepackage{url}
\usepackage{indent first}
\usepackage[labelfont=bf,labelsep=period,justification=raggedright]{caption}
\usepackage[english]{babel}
\usepackage[utf8]{inputenc}
\usepackage{hyperref}
\usepackage[colorinlistoftodos]{todonotes}
\usepackage{tkz-fct}
\usepackage[margin=1in]{geometry} 
\usepackage{tikz}

\DeclareMathOperator{\Fl}{FF}
\DeclareMathOperator{\Fractal}{Fractal}

\newcommand{\eps}{\varepsilon}

\newcommand{\NN}{\mathbb{N}}

\newcommand{\ZZ}{\mathbb{Z}}

\newcommand{\mcG}{\mathcal{G}}

\newcommand{\mcN}{\mathcal{N}}

\newcommand{\mcP}{\mathcal{P}}

\newcommand{\mfG}{\mathfrak{G}}

\newcommand{\msF}{\mathscr{F}}

\theoremstyle{plain}
\newtheorem{thm}{Theorem}
\newtheorem{lemma}[thm]{Lemma}
\newtheorem{cor}[thm]{Corollary}
\newtheorem{conj}[thm]{Conjecture}
\newtheorem{prop}[thm]{Proposition}

\theoremstyle{definition}
\newtheorem{defn}[thm]{Definition}

\theoremstyle{remark}
\newtheorem{rem}[thm]{Remark}
\newtheorem*{ex}{Example}

\numberwithin{equation}{section}
\numberwithin{thm}{section}

\begin{document}

\title{$\mcP$ Play in \textsc{Candy Nim}}

\author[Mani]{Nitya Mani}
\address{(Mani): Stanford University,Department of Mathematics, Stanford, CA 94305}
\email{nityam@stanford.edu}

\author[Nelakanti]{Rajiv Nelakanti}
\address{(Nelakanti, Rubinstein-Salzedo, Tholen): Euler Circle, Palo Alto, CA 94306}
\email{rnelakanti@gmail.com}

\author[Rubinstein-Salzedo]{Simon Rubinstein-Salzedo}
\email{simon@eulercircle.com}

\author[Tholen]{Alex Tholen}
\email{alextholen3.14@gmail.com}

\date{\today}

\maketitle

\begin{abstract}
\textsc{Candy Nim} is a variant of \textsc{Nim} in which both players aim to take the last candy in a game of \textsc{Nim}, with the added simultaneous secondary goal of taking as many candies as possible. We give bounds on the number of candies the first and second players obtain in $3$-pile $\mcP$ positions as well as strategies that are provably optimal for some families of such games. We also show how to construct a game with $N$ candies such that the loser takes the largest possible number of candies and bound the number of candies the winner can take in an arbitrary $\mcP$ position with $N$ total candies.
\end{abstract}

\section{Introduction}

One of the first serious results in the study of combinatorial games was Bouton's solution to the game of \textsc{Nim} in~\cite{Bouton02}.

\begin{defn}\label{d:nim} \textsc{Nim} is a two-player game played with several piles of stones. In a turn, a player removes some number of stones from one pile. The player taking the last stone wins.
\end{defn}

Beyond its historical interest, the game of \textsc{Nim} is interesting because a wide family of games, the so-called \emph{finite normal-play impartial combinatorial games}, can all be reduced to the game of \textsc{Nim}, thanks to the celebrated \emph{Sprague-Grundy theory}, as first described in~\cite{Sprague35} and~\cite{Grundy39}. In this paper, we describe and study a slight modification of the game of \textsc{Nim}, known as \textsc{Candy Nim}, which is interesting in its own right as being a blend of an impartial combinatorial game and a scoring game. While impartial combinatorial games have been widely studied ever since the time of Bouton, the study of scoring games has only recently attracted interest, for instance in~\cite{LNNS15}, ~\cite{larsson2017games}, ~\cite{johnson2014combinatorial}, and~\cite{Stewart13}.

In any \textsc{Nim} game, either the first player to move or the second player to move must have a winning strategy, but not both. We classify the \textsc{Nim} positions based on which player wins with optimal play.

\begin{defn} We call games in which the first player wins with optimal play $\mcN$ positions. Similarly, we call games in which the second player wins with optimal play $\mcP$ positions. We call $\mcN$ and $\mcP$ the \emph{outcome classes} of the associated games.
\end{defn}

\begin{rem}
In the two-player game of \textsc{Nim}, we refer to the losing player as Luca and the winning player as Windsor.
\end{rem}

It is easy to compute the outcome classes and winning strategy for \textsc{Nim} games, based on the following function.

\begin{defn} We define the function $\oplus$ of $a, b \in \ZZ$, called the \textit{nim-sum} as follows: $a\oplus b$ is given by the \textsf{XOR} of $a$ and $b$. (This is given by writing both $a$ and $b$ in binary and adding without carrying.) 
\end{defn}

\begin{ex} Let us compute $9\oplus 5$: \begin{center}\begin{tabular}{cccccc} 9: && 1 & 0 & 0 & 1 \\ 5: &$\oplus$ & 0 & 1 & 0 & 1 \\ \hline 12 && 1 & 1 & 0 & 0 \end{tabular} \end{center} \end{ex}

\begin{thm}[Bouton~\cite{Bouton02}] The \textsc{Nim} game with piles of size $a_1, \ldots , a_p$ is a $\mcP $ position if and only if $n := a_1 \oplus \cdots \oplus a_p = 0$. If $n \neq 0$, winning moves take the total nim-sum to zero.
\end{thm}

\begin{defn} Given a \textsc{Nim} game $G$, with piles of size $a_1,\ldots,a_p$, we define its \emph{Grundy value} $\mcG(G)$ to be $a_1\oplus a_2\oplus\cdots\oplus a_p$. \end{defn}

Because we have an easily computable winning strategy for \textsc{Nim}, the game is boring to play, especially for Luca. In order to make the game more interesting, we change the objective slightly.

\begin{defn}\label{d:candynim} \textsc{Candy Nim} is a two-player combinatorial game with the same setup and game play as \textsc{Nim}. However, in addition to the primary goal of making the last move (as in \textsc{Nim}), players have a secondary goal of collecting as many stones, or \textit{candies} as possible.
\end{defn}

\begin{rem} Note that in \textsc{Candy Nim}, winning always takes priority over collecting candies. No number of candies can fully compensate for the embarrassment of losing the game. An alternative formulation of \textsc{Candy Nim} is that, at the end of the game, the loser (Luca) must give half of her candies to the winner, so that the winner (Windsor) always ends up with more candies than Luca.
\end{rem}

The game of \textsc{Candy Nim} was first introduced by Michael Albert in~\cite{Albert} in an unpublished set of slides based on a talk he gave at the CMS meeting in Halifax in 2004.\footnote{\url{http://www.cs.otago.ac.nz/research/theory/Talks/CandyNim.pdf}} He observed, among other things, that it is not always optimal for Luca to remove candies from the largest pile and also provided some results on values of \textsc{Candy Nim} games (see~\S\ref{sec:prelims} for a definition of the value of a game).

\begin{rem} In this paper, we focus on the $\mcP$ games. This is because Luca has many more options to play with than Windsor. At every turn, in optimal play, Windsor must bring the nim-sum of all of the pile sizes down to zero. This severely limits the options of the winning player. In many of the positions we will study, Windsor will only have a single move available on each turn. On the other hand, Luca loses no matter what her move is, giving room for optimizing her move with respect to the number of candies she collects. Consequently, her turns are more interesting to consider.
\end{rem}

Here, we study certain classes of $\mcP$ positions in \textsc{Candy Nim}. Throughout, we will assume that Windsor is \emph{forced} to play winning moves in the underlying \textsc{Nim} games, so that losing moves are illegal.

After describing the relevant notation in \S\ref{sec:prelims}, we turn to a study of the 3-pile game in \S\ref{sec:3pilelemmas}--\ref{sec:3pile}.\S\ref{sec:3pilelemmas} contains several lemmas that are helpful in the study of 3-pile games and sometimes \textsc{Candy Nim} more generally. In \S\ref{sec:3pilestrats}, we give two strategies for the 3-pile game. The first of these, the \emph{flip-flop strategy}, is a simple strategy for Luca to take as many candies as possible on the current move, subject to allowing Windsor only to remove a single candy. This strategy is easy to work with and analyze, and so is useful for providing bounds for the number of candies each player takes. We also introduce a refinement of the flip-flop strategy, called the \emph{fractal strategy}. The fractal strategy scores better than the flip-flop strategy, and we  conjecture that it is optimal for games in a certain standard form, defined in~\S\ref{sec:prelims}. In~\S\ref{sec:3pile}, we prove bounds on the values of a certain important family of 3-pile games. In Theorem~\ref{t:standard}, we give a strategy for a broad class of $3$-pile games and prove this strategy optimizes the number of candies collected by Luca to within a constant factor of the smallest pile size. This enables us to give fairly strong upper and lower bounds on the difference in candies collected by Luca and Windsor for a large family of $3$-pile games in Theorem~\ref{t:genbound}.

In \S\ref{sec:multipile}, we consider the following problem: Luca gets to distribute an even number $N$ of candies among several piles, subject to the constraint that the resulting position must be a $\mcP$ position. Her goal is to maximize the number of candies she can take.  In Theorem~\ref{thm:equalitycases}, we show that in a game $G$ with $N$ candies, Luca cannot take more than $N - \lfloor \log_2(N) \rfloor$ candies independent of arrangement. In this result, we give explicit characterization of all games where this upper bound is achieved. More generally, we show in Theorem~\ref{thm:5pile} that Luca can always arrange $N$ candies so that Windsor takes at most $O(\sqrt{N})$ candies, in an arrangement with at most $5$ piles. In Theorem~\ref{thm:k-1}, we show that in most games $G$ with $k$ piles, Windsor can take at least $k-1$ candies. 

In \S\ref{sec:examples}, we conclude with some remarks on the 4-pile game, including examples were the best move for Luca is not in the largest pile. We also present some conjectures that we hope will inspire future work.

Throughout this work, we provide some worked examples of optimal play for specific \textsc{Candy Nim} games. We encourage readers to generate their own examples using our program, which is available for download.\footnote{\url{https://github.com/nmani2/candynim}}



\section{Preliminaries} \label{sec:prelims}

We begin with some definitions and notation that will be helpful for the analysis of the $3$ and $n$-pile \textsc{Candy Nim} games. Unless otherwise specified, $G$ will always refer to \textsc{Candy Nim} games.

\begin{defn}
Given a \textsc{Candy Nim} game $G$, let $N(G)$ be the total number of candies in the game. Let $N_W(G)$ be the number of candies collected by winning player Windsor, and let $N_L(G)$ be the number of candies collected by losing player Luca, assuming optimal play.
\end{defn}

Our primary goal is to bound the number of candies Luca (the losing player) can collect relative to Windsor in a $\mcP$ game, assuming optimal play. This difference in candies will be denoted the \textit{value} of the associated \textsc{Candy Nim} game:

\begin{defn} The \textit{value} of a game $V(G)$, is given by \[V(G) = N_L(G) - N_W(G)\]
\end{defn}

\begin{defn}A \textit{turn} is a triple of games $T = (G, G', G'')$ where Luca moves from $G$ to $G'$ and Windsor moves from $G'$ to $G''$. We call each move made by a player from $G$ to $G'$ a \textit{ply} $P = (G, G')$.
\end{defn}

Many of the bounds, such as those in Theorem \ref{t:standard} and Theorem \ref{t:genbound} will arise from the analysis of specific strategies or sequences of moves by Luca and Windsor. To simplify these analyses, we introduce some notation: 

\begin{defn}
The \textit{single-turn value} of a turn $T = (G, G', G'')$, $V_T(G)$ is given as \[V_T(G) = (N(G) - N(G')) - (N(G') - N(G'')) = N(G) + N(G'') - 2N(G').\]
\end{defn}

\begin{defn}
A \textit{strategy} $S$ of game $G$ is a sequence of turns $T_i = (G_i, G_i', G_i'')$ for $1\le i\le n$ such that $G_1 = G$ and for each $j<n$, we have $G_j'' = G_{j+1}$. Furthermore, $G_n''=\varnothing$. We call the \textit{strategic value} of $G$ with strategy $S$ the difference between the number of candies collected by Luca and Windsor under strategy $S$, i.e.\ $V_S(G) = \sum V_{T_i}(G)$. An \textit{optimal strategy} is a strategy such that $\sum V_{T_i}(G) = V(G)$.
\end{defn}

\begin{defn}
The \textit{semiratio} of a turn $T = (G, G', G'')$, $r_T(G)$, is defined to be \[r_T(G) = \frac{N(G) - N(G')}{N(G') - N(G'')}.\]
\end{defn}

\begin{defn} We denote a game $G=[a_1,a_2,a_3,a_4,\ldots,a_p]$ if $G$ is the game with $p$ piles where pile $i$ has $a_i \geq 0$ candies.
\end{defn}

In Theorem \ref{t:standard} and Theorem \ref{t:genbound}, we provide lower and upper bounds on $V(G)$ when $G = [a_1, a_2, a_3]$ is a $3$-pile \textsc{Candy Nim} game. In giving such results, it is helpful to use an alternative characterization of $G$:

\begin{defn}\label{d:g3}
Let $\mfG(a, m, x)$ be the $3$-pile \textsc{Candy Nim} game \[\mfG(a, m, x) = [a,\, 2^{k+1}\cdot m + x,\, 2^{k+1}\cdot m + a \oplus x]\] where $k = \lfloor \log_2 a \rfloor$, $m \geq 1$, and $0 \leq x < 2^k$.
\end{defn}

\begin{rem}
Note that as per Definition~\ref{d:g3}, $2^{k+1}$ is the smallest power of $2$ strictly greater than $a$.
\end{rem}

\begin{defn}
A game $G$ is in \textit{standard form} if \[G = \mfG(2^{k+1} - 1, m, 0) = [2^{k+1} - 1, \,2^{k+1} \cdot m,\, 2^{k+1} (m + 1) - 1].\]
\end{defn}
\begin{defn} Games $G=[g_1,g_2,g_3,\ldots,g_p]$ and $H=[h_1,h_2,h_3,\ldots ,h_q]$ have sum $G+H$ defined by concatenation as follows
\[G+H=[g_1,g_2,g_3,\ldots,g_p,h_1,h_2,h_3,\ldots,h_q] \]
\end{defn}

To gain some basic intuition about \textsc{Candy Nim} consider the following lemma:

\begin{lemma} \label{lem:valuegeq0}
For any game $G$, $V(G) \geq 0$.
\end{lemma}

\begin{proof}
Consider the strategy where Luca takes all of the candies from the largest pile in game $G_i$ for the $i^\text{th}$ turn $T_i$. Then, $V_{T_i}(G_i) \geq 0$ for all $i$, implying that $V(G) \geq 0$.
\end{proof}

\section{Lemmas for the 3-Pile Game} \label{sec:3pilelemmas}

\begin{lemma} \label{lem:oddwinningmoves}
For any $\mcN$ position $G = [g_1, \ldots, g_p]$, suppose there are exactly $j$ piles $g_w$ such that there exists a winning move in pile $w$. Then $j$ is odd.
\end{lemma}

\begin{proof} We give the binary representations of $n$ and $g_i$, as follows:
let $$n = g_1 \oplus \cdots \oplus g_p = \sum _{i = 0} ^k n_i \cdot 2^i, \quad n_k = 1,\, n_i \in \{0,1\}, 1 \le i < k,$$ and for each $1 \leq i \leq p$, let $$g_i = \sum _{h = 0} ^\infty b_{i,h} \cdot 2^h, \quad b_{i,h} \in \{0,1\}.$$ Then Winston has a winning move in pile $g_w$ if and only if $b_{w,k} = 1$. Since $n_k=1$, the number of $w$'s such that $b_{w,k}=1$ must be odd. The piles $g_w$ that contain winning moves are exactly those piles such that $b_{w,k}=1$ (see e.g.~\cite[Proof of Theorem 7.12]{ANW07}), so there are an odd number of piles containing winning moves.
\end{proof}

\begin{lemma} \label{lem:unique3pilemove}
For any 3-pile $\mcP$ position $G$ and a ply $P = (G, G')$ by Luca, there exists a unique $\mcP$ position $G''$ such that $T = (G, G', G'')$ is a turn in \textsc{Candy Nim}. 
\end{lemma}

\begin{proof}
Suppose $G = [g_1, g_2, g_3]$ and Luca moves to $G' = [g_1', g_2, g_3]$. By Lemma~\ref{lem:oddwinningmoves}, Windsor has an odd number of winning moves from $G'$, and he has at most one winning move per pile. Let $n'=\mcG(G')$. Suppose first that Windsor attempts to take $g_1' - g_1''$ candies from the first pile, leaving the game $G'' = [g_1'', g_2, g_3]$. Since $g_2 \oplus g_3 = g_1$ and $g_1 \neq g_1''$, $n'' = g_1'' \oplus g_2 \oplus g_3 = g_1'' \oplus g_1 \neq 0$. Thus Windsor may not move in the first pile. It follows that Windsor may only move in $g_2$ or $g_3$. Since his number of winning moves is odd and at most 2, he has a unique winning move from $G'$.
\end{proof}

\begin{lemma} \label{lem:Semiratio}
Given $G = \mfG(a, m, x)$, and any turn $T = (G, G', G'')$ the semiratio $r_T(G)$ is at most $2a + 1$.
\end{lemma}

\begin{proof}
Consider a turn $T = (G, G', G'')$. We show $r_T(G) \le 2a + 1$. If Luca's ply $(G, G')$ is in the smallest pile, she would take at most $a$ candies, yielding a semiratio at most $a< 2a + 1$.
Suppose Luca moves in either the middle pile or the largest pile such that $G'' = \mfG(a', m', x')$ where $a' \neq a$. There are several cases to consider.
\begin{description}
\item[Case 1] Suppose $G' = [a, 2^{k+1}m + x, a']$ and $G'' = [a, a \oplus a', a']$. Then, \begin{align*} V_T(G) &= (2^{k+1}m + x\oplus a - a') - (2^{k+1}m + x -a'\oplus a) \\ &= x\oplus a - x + a' \oplus a - a' \\ &\leq x+a-x+a'+a-a' \\ &= 2a.\end{align*} Therefore, under this strategy, $$r_T(G) = \frac{N(G) - N(G')}{N(G') - N(G'')} \le \frac{2a+1}{1} \leq 2a + 1.$$
\item[Case 2] Suppose $G' = [a, a', 2^{k+1}m + x\oplus a]$ and $G'' = [a, a \oplus a', a']$. Then, \begin{align*} V_T(G) &= (2^{k+1}m + x - a') - (2^{k+1}m + x\oplus a -a'\oplus a) \\ &= x - a' - x\oplus a + a \oplus a' \\ &\leq x-a'-(a-x)+(a+a') \\ &= 2x.\end{align*} Therefore, $r_T(G) \leq 2x + 1 \le 2a + 1$ under this strategy.
\item[Case 3] Suppose $G' = [a, 2^{k+1}m + x, 2^{k+1}m + x\oplus a']$ and $G'' = [a', 2^{k+1}m + x, 2^{k+1}m + x\oplus a']$. Then, \begin{align*} V_T(G) &= (x\oplus a - x\oplus a') - (a - a') \\ &\leq (x+a)-|x-a'|-a+a' \\ &= a' + x - |a' - x|.\end{align*} This single-turn value is either $2x$ or $2a'$, so $r_T(G) \le 2 \max(x, a') + 1 \leq 2a + 1$ under this strategy.
\item[Case 4] Suppose $G' = [a, 2^{k+1}m + x\oplus a\oplus a', 2^{k+1}m + x\oplus a]$ and $G'' = [a', 2^{k+1}m + x\oplus a\oplus a', 2^{k+1}m + x\oplus a]$. Then, $$V_T(G) = x - x\oplus a\oplus a' - (a - a') \leq x - a + a'.$$ Therefore, $r_T(G) \leq 2a + 1$ under this strategy.
\end{description}
The last possible situation is when $G'' = \mfG(a, m', x')$. Here, there are two cases to consider.
\begin{description}
\item[Case 1] $m' = m$ and $x' < x$. Then, $$N_L(G) - N_L(G'') < 2a + 1.$$
\item[Case 2] $m < m'$. Then, $V_T(G)$ is maximized when $G' = [a, 2^{k+1}m + x, 2^{k+1}m' + x']$. Here, 
\begin{align*}
V_T(G) &= ((2^{k+1}m + x\oplus a) - (2^{k+1}m' + x')) - ((2^{k+1}m + x) - (2^{k+1}m' + x'\oplus a)) \\
&= x\oplus a - x' - x + x'\oplus a \\
&\leq x+a-x'-x+x'+a \\ &= 2a.
\end{align*}
In both cases, $r_T(G) = \frac{N(G) - N(G')}{N(G') - N(G'')} \le 2a+1$. 
\end{description}
\end{proof}

\section{Strategies for the $3$-Pile Game} \label{sec:3pilestrats}

In this section, we present two strategies for certain families of the 3-pile game, which we call the \emph{flip-flop strategy} and the \emph{fractal strategy}. The flip-flop strategy is a simple strategy that, until the last turn, allows Luca to take as many candies as possible subject to allowing Windsor only to take one candy on that turn. The fractal strategy is an iterative variant of the flip-flip strategy which scores better, but for which it is difficult to compute the precise value. We conjecture that some version of the fractal strategy is optimal for games in standard form.

\subsection{The Flip-Flop Strategy} \label{sec:flipflop}

We begin by considering a class of games $[1, 2m, 2m+1]$. We show that they have a simple optimal strategy.
\begin{prop} \label{prop:V1flipflop}
$V([1, 2m, 2m+1]) = 2m$.
\end{prop}
\begin{proof}[Proof] We show optimality and give the strategy inductively,
If $m=1$, then $G = [1,2,3]$. Optimal play occurs when the first move is $T_1 = (G, [1,2,0], [1,1,0])$ with $V(G) = 2$ and $N_L(G) = 4$.
Now assume that $V([1,2m,2m+1])=2m$ is true for all $1 \le m \le m'$. We first show that for $G=[1,2(m'+1),2(m'+1)+1]$, $V(G) \geq 2m'+2$. Consider the strategy $S$ consisting of initial turn $$T_1 = (G, G' = [1,2(m'+1), 2m'], G''= [1,2m' + 1, 2m'])$$ 
followed by optimal play as per the inductive hypothesis for the resulting game $G'' = [1, 2m', 2m' + 1]$. $V_{T_1}(G) = 2$ and $V(G'') = 2m'$ by the inductive hypothesis, giving $V(G) \ge 2m' + 2$. \par 
To show this strategy is optimal, we prove $V(G) \leq 2m'+2$. Consider the four possible cases for Luca's first move.
\begin{description}
\item[Case 1] Consider strategy $S_1$ where Luca takes from the smallest pile. Then the first turn $T_1 = (G, G', G'')$ satisfies $G'' = \{0,2(m'+1),2(m'+1)\}$ by Lemma~\ref{lem:unique3pilemove}. Then, $$V_{S_1}(G) = V_{T_1}(G) + V(G'') = 0 + 0 < 2m' + 2.$$
\item[Case 2] Consider strategy $S_2$ where Luca takes $2k$ candies from the largest pile such that the first turn is $$T_1 = (G, G' = [1,2(m'+1),2j+1], G'' = [1,2j,2j + 1]),$$ where $j = m'+1 - k$. Note that $V_{T_1}(G) = 0$. By induction, $V(G'') = 2j$. Therefore, $$V_{S_2}(G) = V_{T_1}(G) + V(G'') = 2j < 2m'+2.$$
\item[Case 3] Consider strategy $S_3$ where Luca takes $2k + 1$ candies from the largest pile such that the first turn is $$T_1 = (G, G' = [1,2(m'+1),2j], G'' = [1,2j+1,2j]),$$ where $j = m'+1 - k$. This time, $V_{T_1}(G) = 2$. By induction, $V(G'') = 2j$. Therefore, $$V_{S_3}(G) = V_{T_1}(G) + V(G'') = 2 + 2j \leq 2m'+2.$$
\item[Case 4] Consider strategy $S_4$ where Luca takes $k$ candies from the medium pile such that the first turn is $$T_1 = (G, G' = [1, 2j, 2(m'+1) + 1], G'' = [1, 2j, 2j+1]$$ or $$T_1 = (G, G' = [1, 2j + 1, 2(m'+1) + 1], G'' = [1, 2j, 2j+1].$$ In this case, $V_{T_1} \leq 0$. By induction, $V(G'') = 2j$. Therefore, $$V_{S_4}(G) = V_{T_1}(G) + V(G'') \leq 2j < 2m'+2.$$
\end{description}
Since $\max(V_{S_1}(G), V_{S_2}(G), V_{S_3}(G), V_{S_4}(G) \le 2m'+2$, we obtain the desired equality $V(G) = 2m' +2$. \end{proof}

\begin{figure} [h]
\begin{displaymath}
    \xymatrix{ [1 & 2m  & 2m+1\ar@[red][d]^{\,\color{red}(-3)}]\\
               [1 & 2m\ar@[green][d]^{\color{green}\,(-1)} & 2(m-1)] \\
               [1 & 2(m-1)+ 1\ar@[red][d]^{\color{red}\,(-3)} & 2(m-1)] \\
               [1 & \vdots & \vdots\ar@[green][d]^{\color{green}\,(-1)}] \\
               [1 & 2 & 3\ar@[red][d]^{\color{red}\,(-3)}] \\ [1 & 2\ar@[green][d]^{\color{green}\,(-1)} & 0] \\ [1 & 1 & 0] }
\end{displaymath}
\caption{The sequence of moves that occurs when using the flip-flop strategy in the game $G = [1, 2m, 2m+1]$.}
\label{fig:flipflop}
\end{figure}
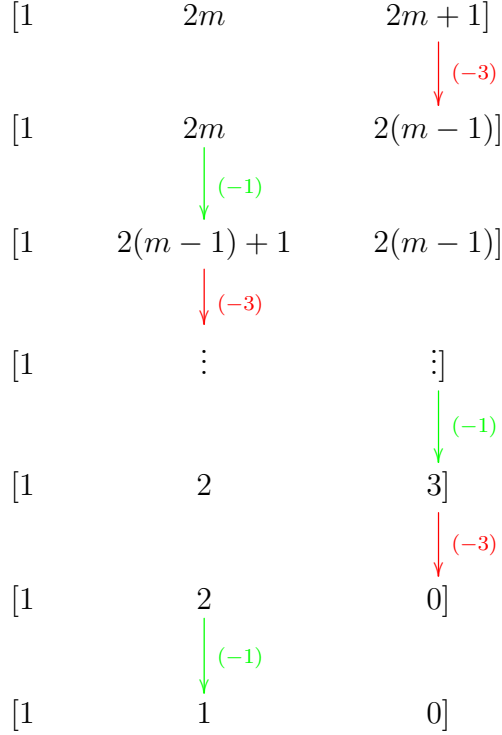

The inductive optimal strategy in the proof above is pictorially represented by the sequence of moves in Figure~\ref{fig:flipflop}.

More generally, we can consider the following strategy, for games of the form $G = [2^k-1,2^k \cdot m,2^k \cdot (m+1)-1]$, that we will term the \textit{flip-flop} strategy:

\begin{defn}\label{d:flip-flop}
Given a game $G = \mfG(2^k-1,m,0)=[2^k-1,2^k \cdot m,2^k \cdot (m+1)-1]$, the \textit{flip-flop strategy} $\Fl(G)$ is given as follows: \begin{enumerate} \item If $m\ge 1$, Luca removes $2^{k+1}-1$ candies from the third pile, then Windsor removes one candy from the second pile. The resulting game is $\mfG(2^k-1,m-1,0)$. Then continue with $\Fl(\mfG(2^k-1,m-1,0))$. \item If $m=0$, then we have $G=(2^k-1,2^k-1)$. Luca removes one pile, then Windsor removes the other one. \end{enumerate}
\end{defn}

\begin{prop} \label{prop:4.3}
For $G = [2^k-1,2^k \cdot m,2^k \cdot (m+1)-1]$, we have
$$V(G) \ge V_{\Fl(G)}(G) = (m-1) \cdot (2^{k+1} - 2).$$
\end{prop}

\begin{proof}
Consider the strategy $\Fl(G)$ with initial turn $T_1 = (G = G^{(0)}, G^{(1)}, G^{(2)})$, where Luca takes $2^{k+1} - 1$ candies from the largest pile and Windsor takes $1$ candy from the middle pile. Thus, $V_{T_0}(G) = 2^{k+1} - 2$ with $G^{(2)} = [2^k - 1, 2^k \cdot (m-1), 2^k \cdot m - 1]$. Repeat for turns $T_2, \ldots, T_{m-1}$ where $$T_i = (G^{(2i-2)},G^{(2i-1)} ,G^{(2i)}).$$
For $i = 1, \ldots, m-1$, $$V_{T_i}(G^{(2i-2)}) = 2^{k+1} - 2.$$
When $m = 0$, the resulting game is $G^{(2m-2)} = [2^k - 1, 2^k - 1]$ with $V(G^{(2m-2)}) = 0$. Thus, $V_{\Fl(G)}(G) = (m-1) \cdot (2^{k+1} - 2) + 0$.
\end{proof}

\subsection{The Fractal Strategy} 
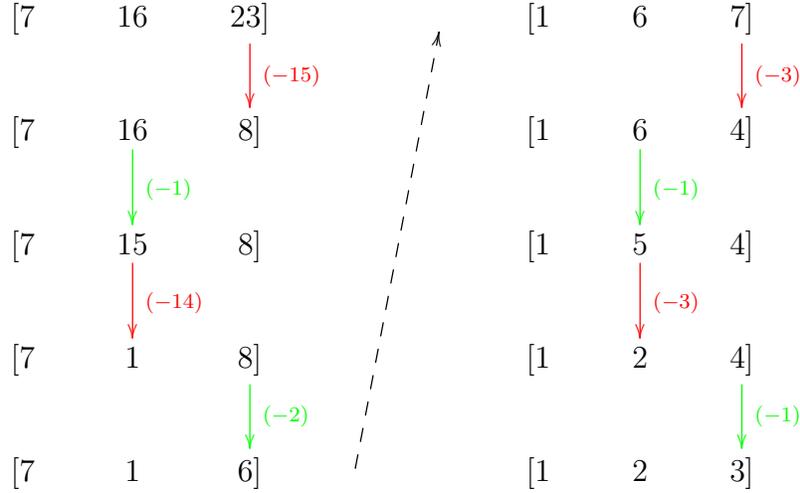
\begin{figure}
\begin{displaymath}
    \xymatrix{ [7 & 16  & 23\ar@[red][d]^{\,\color{red}(-15)}] &\,&\,& [1 & 6 & 7\ar@[red][d]^{\,\color{red}(-3)}]\\
               [7 & 16\ar@[green][d]^{\color{green}\,(-1)} & 8] &\,&\,& [1 & 6\ar@[green][d]^{\color{green}\,(-1)} & 4]\\
               [7 & 15\ar@[red][d]^{\color{red}\,(-14)} & 8] &\,&\,& [1 & 5\ar@[red][d]^{\,\color{red}(-3)} & 4]\\
               [7 & 1 & 8\ar@[green][d]^{\color{green}\,(-2)}] &\,&\,& [1 & 2 & 4\ar@[green][d]^{\color{green}\,(-1)}]\\
               [7 & 1 & 6] &\ar@{-->}[ruuuu]\,&\,& [1 & 2 & 3] }
\end{displaymath}
\caption{In the fractal strategy with $G=[7,16,23]$, we begin by applying the flip-flop strategy until the game reaches $H = [7, 15, 8]$. If Luca continued via the flip-flop strategy, the next turn would be $T = (H, [7, 8], [7,7])$, giving Luca $22/30$ candies in $H$. If Luca instead reduced the smallest pile from size $7$ to $1$, the single-turn value of that reduction would be $12$. This yields the game $[1, 6, 7]$ which has value $6$, as shown in Proposition~\ref{prop:V1flipflop}. With this sequence of moves, Luca does better, obtaining $24/30$ candies of $H$}
\label{fig:fractal1}
\end{figure}

We can improve the above strategy to one that exhibits a curious fractal-like behavior, as in Figure~\ref{fig:fractal1} and the following example:

\begin{prop}
$62m + 60 \geq V([31, 32m, 32m + 31]) \geq 62(m-1) + 98$
\end{prop}

\begin{proof}\
\begin{description}
\item[Upper Bound] By Lemma \ref{lem:Semiratio}, the semiratio on any turn is at most $63$, implying that $$V([31, 32m, 32m + 31]) \le \frac{63-1}{63+1}\cdot (31 + 32m + 32m + 31) = 60.0625 + 62m.$$

\item[Lower Bound] Consider the following strategy, broken down into $m>1$ and $m=1$.
\begin{enumerate}
\item While $m>1$, Luca recursively removes $63$ candies from the largest pile, requiring Windsor to respond by removing $1$ candy from the middle pile, creating the turn $$T = (G, [31, 32m, 32(m-1)], [31,32(m-1) + 31, 32(m-1)).$$ The accumulated value is $\sum_{m = 2}^{m} (63-1)=62(m-1)$.
\item When $m = 1$, the turn is $$T_4 = ([31,32,63],[31,32,3],[31,28,3])$$ with a single-turn value of $V_{T_4}(G_4) = 60-4 = 56$. \begin{enumerate} \item Then, for all games $G_2 = [3, 4m', 4m' + 3]$ with $m' > 1$, the turn would be $$T = (G_2, [3,4m', 4(m'-1)], [3, 4(m'-1) + 3, 4(m'-1)]).$$ The accumulated value for all $G_2$ is $\sum _{m'>1} (7-1) = 6(7-1) = 36$. \item When $m'=1$, the turn is $T_1 = ([3,4,7],[3,4,1],[3,2,1])$. Finally, the last two turns are $T_0 = ([1,2,3],[1,2],[1,1])$ and $T_0' = ([1,1],[1],\varnothing)$. In total, we get an overall value of $62(m-1) + 56 + 36 + 4 + 2 + 0 = 62(m-1) + 98$. \end{enumerate} \end{enumerate} \end{description}
\end{proof}

\begin{rem}
By checking numerically, we have verified that for $m < 12$, $V([31, 32m, 32m + 31]) = 62(m-1) + 98$. We conjecture equality holds for all $m \in \NN$. 
\end{rem}

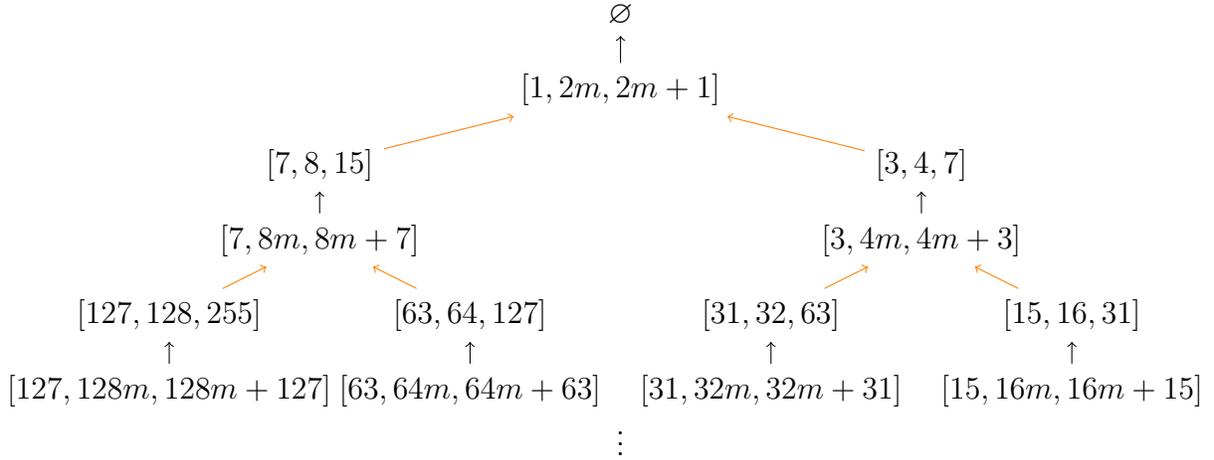
\begin{figure}[h]
\begin{tikzpicture}
	\node (0) at (0,8) {$\varnothing$};
    \node (05) at (0,7) {$[1,2m,2m+1]$};
    \draw [->] (05) -- (0);
    
	\node (1) at (4,6) {$[3,4,7]$};
	\node (2) at (-4,6) {$[7,8,15]$};
    \draw [orange] [->] (1) -- (05);
    \draw [orange] [->] (2) -- (05);
    
    \node (15) at (4,5) {$[3, 4m, 4m + 3]$};
    \node (25) at (-4,5) {$[7, 8m, 8m + 7]$};
    \draw [->] (15) -- (1);
    \draw [->] (25) -- (2);
    
    \node (3) at (6,4) {$[15,16,31]$};
    \node (4) at (2,4) {$[31,32,63]$};
    \draw [orange] [->] (3) -- (15);
    \draw [orange] [->] (4) -- (15);
    
    \node (5) at (-2,4) {$[63,64,127]$};
    \node (6) at (-6,4) {$[127,128,255]$};
  	\draw [orange] [->] (5) -- (25);
    \draw [orange] [->] (6) -- (25);
    
    \node (35) at (6,3) {$[15,16m,16m+15]$};
    \node (45) at (2,3) {$[31,32m,32m+31]$};
    \node (55) at (-2,3) {$[63,64m,64m+63]$};
    \node (65) at (-6,3) {$[127,128m,128m+127]$};
    
    \draw [->] (35) -- (3);
    \draw [->] (45) -- (4);
    \draw [->] (55) -- (5);
    \draw [->] (65) -- (6);

	\node (d) at (0,2.4) {$\vdots$};
    
\end{tikzpicture}
\caption{We represent the fractal strategy using black and orange arrows. The black arrows indicate the flip-flop strategy of Section~\ref{sec:flipflop} and the orange arrows indicate a change of smallest pile size.}
\end{figure}

\begin{defn}
Define a function $f: \NN \rightarrow \NN$ to be \textit{contractive} if for all $a \in \NN$, $f(a) \le a$. Let $\msF$ denote the family of contractive functions.
\end{defn}

\begin{defn} \label{d:fractal}
Consider a game $G$ of the form $G = \mfG(2^k-1,m,0)=[2^k - 1, 2^k \cdot m, 2^k (m+1) -1]$, for $m, k \ge 1$. Let $f \in \msF$. We define the \textit{fractal strategy} $\Fractal_f(G)$ based on $f$ as follows:
\begin{enumerate}
\item If $m>1$, then Luca plays as in $\Fl$ by removing $2^{k+1}-1$ candies from the third pile, and then Windsor moves to $\mfG(2^k-1,m-1,0)$. Then play $\Fractal_f(\mfG(2^k-1,m-1,0))$.
\item If $m=1$ and $f(a)=a$, then play as in the flip-flop strategy.
\item If $m=1$ and $f(a)<a$, then Luca moves the smallest pile to $2^{f(a)}-1$, and Windsor moves to $\mfG(2^{f(a)}-1,2^{a-f(a)},0)$. Then play $\Fractal_f$ from there. 
\end{enumerate}
\end{defn}

\begin{thm} 
Given $G = \mfG(2^k - 1, m, 0)$ with $k, m \ge 1$, we have
$$\sup_{f \in \msF} V_{\Fractal_f}(G) = (m-2) \cdot (2^{k+1} - 2) + \displaystyle \sum_{i=0}^{\lceil \log_2 k \rceil } 2^{\lfloor \frac k{2^i}\rfloor+1}-2^{\left \lfloor \frac k{2^{i+1}} \right \rfloor+1} + \left (2^{\left \lfloor \frac k{2^{i+1}}\right \rfloor +1}-1 \right ) \left ( 2^{\left \lfloor \frac k{2^i}\right \rfloor - \left \lfloor \frac k{2^{i+1}}\right \rfloor}-2\right ),$$
with the supremum achieved by taking $f: a \mapsto \lfloor\frac{a}{2}\rfloor$.

\end{thm}
\begin{proof}
We first show that the fractal strategy $f:a\mapsto\lfloor\frac{a}{2}\rfloor$ achieves the stated bound. Consider the strategy $\Fractal_f$ where $f: a \mapsto \lfloor \frac a2 \rfloor$. We show \[V_{\Fractal_f}(G) = (m-2) \cdot (2^{k+1} - 2) + \displaystyle \sum_{i=0}^{\lceil \log_2 k \rceil } 2^{\lfloor \frac k{2^i}\rfloor+1}-2^{\left \lfloor \frac k{2^{i+1}} \right \rfloor+1} + \left (2^{\left \lfloor \frac k{2^{i+1}}\right \rfloor +1}-1 \right ) \left ( 2^{\left \lfloor \frac k{2^i}\right \rfloor - \left \lfloor \frac k{2^{i+1}}\right \rfloor}-2\right ).\]
It suffices to show that for the game $H=[2^k-1,2^k,2^{k+1}-1]$, 
\[V_{\Fractal_f}(H) = \displaystyle \sum_{i=0}^{\lceil \log_2 k \rceil } 2^{\lfloor \frac k{2^i}\rfloor+1}-2^{\left \lfloor \frac k{2^{i+1}} \right \rfloor+1} + \left (2^{\left \lfloor \frac k{2^{i+1}}\right \rfloor +1}-1 \right ) \left ( 2^{\left \lfloor \frac k{2^i}\right \rfloor - \left \lfloor \frac k{2^{i+1}}\right \rfloor}-2\right ).\]
Under $\Fractal_f$, the first turn is \[T = (H, H', H'') = ([2^k-1,2^k,2^{k+1}-1], [2^k-1,2^k,2^{\lfloor \frac k2 \rfloor}-1], [2^k-1,2^k-2^{\lfloor \frac k2 \rfloor},2^{\lfloor \frac k2 \rfloor}-1]).\] 
Under the this strategy, we perform $\Fl$ until we reach the game $ [2^{\lfloor \frac k2 \rfloor}-1,2^{\lfloor \frac k2 \rfloor},2^{\lfloor \frac k2 \rfloor-1}]$. This involves repeating the following sequence of moves $2^{k-\lfloor \frac k2 \rfloor}-2 $ times: $$[2^{\lfloor \frac k2 \rfloor}-1,a \cdot 2^{\lfloor \frac k2 \rfloor}, (a+1) \cdot 2^{\lfloor \frac k2 \rfloor} -1] \mapsto [2^{\lfloor \frac k2 \rfloor}-1,a \cdot 2^{\lfloor \frac k2 \rfloor}, (a-1) \cdot 2^{\lfloor \frac k2 \rfloor} ] \mapsto [2^{\lfloor \frac k2 \rfloor}-1,a \cdot 2^{\lfloor \frac k2 \rfloor}-1, (a-1) \cdot 2^{\lfloor \frac k2 \rfloor} ].$$
Since for any fractal strategy $g$ when $a \neq 1$ $$V_{\Fractal_g}([2^{\lfloor \frac k2 \rfloor}-1,a \cdot 2^{\lfloor \frac k2 \rfloor}, (a+1) \cdot 2^{\lfloor \frac k2 \rfloor} -1])=V_{\Fractal_g}([2^{\lfloor \frac k2 \rfloor}-1,a \cdot 2^{\lfloor \frac k2 \rfloor}-1, (a-1) \cdot 2^{\lfloor \frac k2 \rfloor} ])+2^{\lfloor \frac k2 \rfloor +1}-1,$$ we obtain that
$$V_{\Fractal_f}(H)=V_{\Fractal_f}([2^{\lfloor \frac k 2 \rfloor}-1,2^{\lfloor \frac k2 \rfloor},2^{\lfloor \frac k2 \rfloor +1}-1])+2^{k+1}-2^{\lfloor \frac k2 \rfloor}+(2^{\lfloor \frac k2 \rfloor +1}-1)(2^{k-\lfloor \frac k2 \rfloor}-2).$$
Via the inductive hypothesis we obtain the desired result: \[V_{\Fractal_f}(G) = (m-2) \cdot (2^{k+1} - 2) + \displaystyle \sum_{i=0}^{\lceil \log_2 k \rceil } 2^{\lfloor \frac k{2^i}\rfloor+1}-2^{\left \lfloor \frac k{2^{i+1}} \right \rfloor+1} + \left (2^{\left \lfloor \frac k{2^{i+1}}\right \rfloor +1}-1 \right ) \left ( 2^{\left \lfloor \frac k{2^i}\right \rfloor - \left \lfloor \frac k{2^{i+1}}\right \rfloor}-2\right ).\]

We now show that the $\Fractal_f$ strategy for $f: a \mapsto \lfloor \frac a2 \rfloor$ is optimal over all possible strategies $\Fractal_g$. It suffices to show that for all fractal strategies $\Fractal_g$, $$V_{\Fractal_g}([2^k-1,2^k,2^{k+1}-1]) \le V_{\Fractal_f}([2^k-1,2^k,2^{k+1}-1]).$$ We consider two cases based on whether $g(k) < \lfloor \frac k2 \rfloor$ or $g(k) > \lfloor \frac k2 \rfloor$ 
We first show that if $g(k)=i<\lfloor \frac k2 \rfloor$, then
\begin{equation} V_{\Fractal_g}([2^k-1,2^k,2^{k+1}-1]) \le V_{\Fractal_f}([2^k-1,2^k,2^{k+1}]). \label{eq:gvsffractal1}\end{equation}

  Equivalently, we wish to show that the left side of~(\ref{eq:gvsffractal1} minus the right side is less than or equal to zero. 
By the definition of the fractal strategy, we have \begin{align*}
&V_{\Fractal_g}([2^k-1,2^k,2^{k+1}-1])-V_{\Fractal_f}([2^k-1,2^k,2^{k+1}]) \\
&\le V_{\Fractal_g}([2^i-1,2^i,2^{i+1}-1]) - V_{\Fractal_f}([2^{\lfloor \frac k2 \rfloor}-1,2^{\lfloor \frac k2 \rfloor},2^{\lfloor \frac k2 \rfloor+1}-1])\\
&\qquad \qquad -2^{k-i} + 2^{i+1} +2^{k - \lfloor \frac k2 \rfloor} +2^{\lfloor \frac k2 \rfloor -i}-2\\ 
&\le V_{\Fractal_f}([2^i-1,2^i,2^{i+1}-1])-V_{\Fractal_f}([2^i-1,2^i,2^{i+1}-1]) \\
&\qquad \qquad + 2^{i+1} + 2^{k-\lfloor \frac k2 \rfloor} + 2^{\lfloor \frac k2 \rfloor -i} -2^{k-i}-2\\ 
&=2^{i+1}+2^{k-\lfloor \frac k2 \rfloor} + 2^{\lfloor \frac k2 \rfloor -i} - 2^{k-i}-2. \end{align*}
Now, suppose $i \neq {\lfloor \frac k2 \rfloor}-1$. Then we have \begin{align*}
&V_{\Fractal_g}([2^k-1,2^k,2^{k+1}-1])-V_{\Fractal_f}([2^k-1,2^k,2^{k+1}-1]) \\
&\le 2^{i+1}+2^{k-\lfloor \frac k2 \rfloor} + 2^{\lfloor \frac k2 \rfloor -i} - 2^{k-i} -2 \\ &\le 2^{\lfloor \frac k2 \rfloor-1}+2^{\lfloor \frac k2 \rfloor+1} +2^{\lfloor \frac k2 \rfloor} - 2^{k-\lfloor \frac k2 \rfloor +2} -2\\
&\le 2^{\lfloor \frac k2 \rfloor +2} - 2^{\lfloor \frac k2 \rfloor +2}-2 \\ &< 0.
\end{align*}
On the other hand, if $i=\lfloor \frac k2 \rfloor -1$, we get \begin{align*}
&V_{\Fractal_g}([2^k-1,2^k,2^{k+1}-1])-V_{\Fractal_f}([2^k-1,2^k,2^{k+1}-1]) \\
&\le 2^{i+1}+2^{k-\lfloor \frac k2 \rfloor} + 2^{\lfloor \frac k2 \rfloor -i} - 2^{k-i} -2  \\ &\le 2^{\lfloor \frac k2 \rfloor} +2^{k-\lfloor \frac k2 \rfloor} + 2^{1} - 2^{k-\lfloor \frac k2 \rfloor +1} -2  \\ 
&\le 2^{\lfloor \frac k2 \rfloor} - 2^{k-\lfloor \frac k2 \rfloor} \\ &\le 0.
\end{align*}
Next suppose that
$g(k) = i >\lfloor \frac k2 \rfloor$. First, recall the notation $\mfG(a, m, x) = [a, 2^k\cdot m + x,\, 2^k\cdot m + a \oplus x], 2^k>a \ge 2^{k-1}$. Let $f',g'\in\msF$ be defined by $f'(n)=f(n)$ if $n\neq\lfloor\frac{k}{2}\rfloor$ and $f'(\lfloor\frac{k}{2}\rfloor)=\lfloor\frac{i}{2}\rfloor$, and $g'(n)=f(n)$ if $n\neq k$ and $g'(k)=i$. We will show that $V_{\Fractal_{f'}}(\mfG(2^k-1,1,0))\ge V_{\Fractal_{g'}}(\mfG(2^k-1,1,0))$. By induction, this implies that $V_{\Fractal_{f}}(\mfG(2^k-1,1,0))\ge V_{\Fractal_{g}}(\mfG(2^k-1,1,0))$. To this end, we have 
\begin{align*}
&V_{\Fractal_{f'}}(\mfG(2^k-1,1,0))-V_{\Fractal_{g'}}(\mfG(2^k-1,1,0)) \\
&=V_{\Fractal_{f'}}(\mfG(2^{\lfloor \frac k2 \rfloor}-1,1,0))-V_{\Fractal_{g'}}(\mfG(2^i-1,1,0))-2^{\lfloor \frac k2 \rfloor} \\
&\qquad \qquad +(2^{\lfloor \frac k2 \rfloor +1}-1)(2^{k-\lfloor \frac k2\rfloor}-2)+2^{i}-(2^{i+1}-1)(2^{k-i}-2) \\ 
&=V_{\Fractal_{f'}}(\mfG(2^{\lfloor \frac k2 \rfloor}-1,1,0))-V_{\Fractal_{g'}}(\mfG(2^i-1,1,0)) \\
&\qquad \qquad -2^{\lfloor \frac k2 \rfloor}-2^{k-\lfloor \frac k2 \rfloor} - 2^{\lfloor \frac k2 \rfloor +2}+2^{i} + 2^{k-i} +2^{i+2}  \\
&\ge V_{\Fractal_{f'}}(\mfG(2^{\lfloor \frac k2 \rfloor}-1,1,0)) \\
&\qquad \qquad - V_{\Fractal_{f'}}(\mfG(2^i-1,1,0))-2^{\lfloor \frac k2 \rfloor}-2^{k-\lfloor \frac k2 \rfloor} - 2^{\lfloor \frac k2 \rfloor +2}+2^{i} + 2^{k-i} +2^{i+2} 
\\
&=-2^{k-\lfloor \frac k2 \rfloor} - 2^{\lfloor \frac k2 \rfloor +2} + 2^{k-i} +2^{i+2}+ (2^{\lfloor \frac i2 \rfloor+1}-1)(2^{\lfloor \frac k2 \rfloor-\lfloor \frac i2\rfloor}-2)-(2^{\lfloor \frac i2 \rfloor+1}-1)(2^{i-\lfloor \frac i2 \rfloor}-2) \\
&=-2^{k-\lfloor \frac k2 \rfloor} - 2^{\lfloor \frac k2 \rfloor +2} + 2^{k-i} +2^{i+2}+ 2^{\lfloor \frac k2 \rfloor+1} - 2^{\lfloor \frac k2 \rfloor - \lfloor \frac i2 \rfloor} - 2^{\lfloor \frac i2 \rfloor +2}  - 2^i +2^{\lfloor \frac i2 \rfloor+2} + 2^{i-\lfloor \frac i2 \rfloor} \\ &=2^{k-i}+3 \cdot 2^{i}+ 2^{i-\lfloor \frac i2 \rfloor}- 2^{\lfloor \frac k2 \rfloor - \lfloor \frac i2 \rfloor}-2^{\lfloor \frac k2 \rfloor+1} -2^{k-\lfloor \frac k2 \rfloor}.
\end{align*}
Now, we can see that the trio of inequalities $i \ge k-\lfloor \frac k2 \rfloor$ and $i \ge {\lfloor \frac k2 \rfloor +1}$ and $i \ge \lfloor \frac k2 \rfloor - \lfloor \frac i2 \rfloor$ are each true, and so that allows us to simplify to get\begin{align*}
&V_{\Fractal_f}(\mfG(2^k-1,1,0))-V_{\Fractal_g}(\mfG(2^k-1,1,0))\\
&\ge2^{k-i}+3 \cdot 2^{i}+ 2^{i-\lfloor \frac i2 \rfloor}- 2^{\lfloor \frac k2 \rfloor - \lfloor \frac i2 \rfloor}-2^{\lfloor \frac k2 \rfloor+1} -2^{k-\lfloor \frac k2 \rfloor} \\ &\ge 2^{k-i} + 2^{i-\lfloor \frac i2 \rfloor}, \end{align*} which is positive.
This resolves the last case, yielding the desired result.
\end{proof}

\section{Bounds for the 3-Pile Game} \label{sec:3pile}

\begin{thm} \label{t:standard}
Given a standard form game $G = \mfG(2^{k+1} - 1, m, 0)$, we have \[V(\mfG(2^{k+1} - 1, m, 0)) \le (2^{k+2} - 2)m + (2^{k+2} - 2) - 2 + \delta_{0k},\] where $\delta_{0k}$ is the Kronecker delta function which is 1 if $k=0$ and 0 otherwise. Furthermore, 
\[
V(\mfG(2^{k+1} - 1, m, 0)) \ge 2(2^{k+1} - 1)m - 2(2^{ \lceil \frac{k}{2}\rceil}  - 1)  + V(\mfG(2^{ \lceil\frac{k}{2} \rceil } - 1, 2^{ \lfloor \frac{k}{2} \rfloor + 1} - 1, 0))
\]
Alternatively,
\[
V(\mfG(2^{k+1} - 1, m, 0)) \geq 2(2^{k+1} - 1)(m-1) + b(k),
\] where
\[
3(2^{k+1} - 1) \leq b(k) \leq 4(2^{k+1} - 1) - 2.
\]

\end{thm}

\begin{proof} 
($\le$) We will first show that \[V(G) = V(\mfG(2^{k+1} - 1, m, 0)) \leq (2^{k+1} - 2)m + (2^{k+1} - 2) - 2 + \delta_{0k}.\] By Lemma \ref{lem:Semiratio}, $r_T(G) \le s = 2^{k+2} - 1$. Then,
\begin{align*}
V_L(G) &\leq \frac{s-1}{s+1} N(G) \\
&= \frac{2^{k+2} - 2}{2^{k+2}} \cdot (2^{k+2} - 2 + 2^{k+2}m) \\
&=2^{k+2} - 2 + (2^{k+2} - 2)m-2\frac{2^{k+2} - 2}{2^{k+2}} \\
&= (2^{k+2} - 2)m + (2^{k+2} - 2) - 2 + \frac{1}{2^k}.
\end{align*} 
($\ge$)
Given the game $G = \mfG(2^{k+1} - 1, m, 0)$, consider the strategy where Luca removes $2^{k+2}-1$ candies from the largest pile when $m > 1$ and $2^{k+2} - 2^{\lfloor \frac{k+1}{2} \rfloor}$ from the largest pile when $m = 1$. Then, 
$$V_L(G) \geq 2(2^{k+1} - 1)m - 2(2^{\lceil \frac{k}{2} \rceil} - 1) + V(\mfG(2^{\lceil \frac{k}{2} \rceil} - 1, 2^{\lfloor \frac{k}{2} \rfloor + 1} - 1, 0)).$$ 
This is an example of the fractal strategy as in Definition~\ref{d:fractal} with $f(k) = \lfloor \frac{k+1}{2} \rfloor$. 
\end{proof}

\begin{cor} \label{cor:UsefulBound}
$V(\mfG(a, m, x))\ge 2a(m-1) + x\oplus a + a - x$.
\end{cor}

\begin{proof}
Let $$G = \mfG(a, m, x) = [a, 2^{\lfloor \log_2 a \rfloor + 1}m + x, 2^{\lfloor \log_2 a \rfloor + 1}m + x \oplus a].$$ Consider first turn $T_0 = (G, G', G'')$ such that $$G' = [a, 2^{\lfloor \log_2 a \rfloor + 1}m + x, 2^{\lfloor \log_2 a \rfloor + 1}(m-1)]$$ and $$G'' = [a, 2^{\lfloor \log_2 a \rfloor + 1}(m-1) + a, 2^{\lfloor \log_2 a \rfloor + 1}(m-1)].$$ Then, $V_{T_0}(G) = x\oplus a + a - x$. 
For $0 < i < m$, let the $i^\text{th}$ turn be $T_i = (G_i, G_i', G_i'')$ where  (similar to the flip-flop strategy of Definition~\ref{d:flip-flop}),
$$G_i = [a, 2^{\lfloor \log_2 a \rfloor + 1}(m-i), 2^{\lfloor \log_2 a \rfloor + 1}(m-i) + a],$$
$$G_i' = [a, 2^{\lfloor \log_2 a \rfloor + 1}(m-i), 2^{\lfloor \log_2 a \rfloor + 1}(m- i - 1)],$$
$$G_i'' = [a, 2^{\lfloor \log_2 a \rfloor + 1}(m-i-1) + a, 2^{\lfloor \log_2 a \rfloor + 1}(m-i-1)],$$
with $V_{T_i}(G) = 2a$. 
After turn $T_{m - 1}$, $G_{m-1}'' = [a, a]$. Therefore the game concludes after $m$ turns and in total, Luca takes $2a(m-1) + x\oplus a + a - x$ candies.
\end{proof}

\begin{thm} \label{t:genbound}
If $G = \mfG(2^{k+1} - 1, m, x)$, then 
\[V(\mfG(2^{k+1} - 1, m - 1, 0)) + 2(2^{k+1} - 1) - 2x \leq V(G) \leq V(\mfG(2^{k+1} - 1, m + 1, 0)) - 2(2^{k+1} - 1) + 2x.\]

\end{thm}

\begin{proof} 
Let $G = [2^{k+1} - 1, 2^{k+1}m + x, 2^{k+1}m + 2^{k+1} - 1 - x]$. \par

$(\leq )$. We construct a strategy $S$ that achieves a value of $$V_S(G) = V(\mfG(2^{k+1} - 1, m-1, 0)) + 2(2^{k+1} - 1) - 2x.$$ Let the first turn $T_1 = (G, G', G'')$ consist of $$G' =  [2^{k+1} - 1, 2^{k+1}m + x, 2^{k+1}(m-1)]$$ and $$G'' = [2^{k+1} - 1, 2^{k+1}(m-1) + 2^{k+1} - 1, 2^{k+1}(m-1)].$$
Then the single-turn value is $V_{T_1}(G) = 2(2^{k+1} - 1) - 2x$ and $G'' = \mfG(2^{k+1} - 1, m - 1, 0)$, yielding an overall value of $$V_S(G) = V(\mfG(2^{k+1} - 1, m-1, 0)) + 2(2^{k+1} - 1) - 2x,$$ which gives a lower bound on $V(G)$. \par

$(\geq)$. We prove $$V(\mfG(2^{k+1} - 1, m + 1, 0)) - 2(2^{k+1} - 1) + 2x \geq V(G).$$
Given game $G_0 = \mfG(2^{k+1} - 1, m + 1, 0)$, under any strategy $S$, $$V_S(G_0) \leq V(G_0).$$ 
Suppose the first turn in $S$ is $T_1 = (G_0, G_0', G_0'')$, where $$G_0' = [2^{k+1} - 1, 2^{k+1}(m+1), 2^{k+1}m + x]$$ and $$G_0'' = [2^{k+1} - 1, 2^{k+1}m + x, 2^{k+1}m + 2^{k+1} - 1 - x].$$ Note that $G_0''$ is the only move to a $\mcP$ position from $G_0'$ . The resulting game is $G_0'' = G$, implying that $V_S(\mfG(2^{k+1} - 1, m + 1, 0)) = 2(2^{k+1} - 1) - 2x + V(G)$. This gives the desired inequality: $$ 2(2^{k+1} - 1) - 2x + V(G) = V_S(\mfG(2^{k+1} - 1, m + 1, 0)) \leq V(\mfG(2^{k+1} - 1, m + 1, 0)).$$
\end{proof}

\section{Optimal Allocation of $N$ Candies} \label{sec:multipile}

So far, we have only considered \textsc{Candy Nim} positions with three piles. We have seen that, in these games, Luca can take a substantial majority of the candies, and indeed there are 3-pile $\mcP$ positions in which Luca takes a proportion of at least $1-\eps$ of the candies, for any fixed $\eps>0$. It is natural, then, to consider the problem of Luca allocating $N$ candies, in a $\mcP$ position, so that she maximizes the number of candies that she can take with optimal play. This problem is the motivating question for this section.

\begin{lemma}\label{l:p2} If $G \in \mcP$, then $N(G)-N(G') \le \frac{N(G)}2$.
\end{lemma}

\begin{proof} Let $G=[a_1,a_2,\ldots,a_p]$, where $a_1 \ge a_2 \ge \cdots \ge a_p$. For any ply $(G,G')$, we have $N(G)-N(G')\le a_1$. So, it suffices to show that $a_1  \le \frac {N(G)} 2$. Since $G \in \mcP$, we have $a_1 = a_2 \oplus a_3 \oplus \cdots \oplus a_p$. For any $x_1,\ldots,x_k$, we have $x_1\oplus\cdots\oplus x_k\le x_1+\cdots+x_k$, so $a_1 \le a_2 + a_3 + \cdots + a_p$. Since $N(G)=a_1 + a_2 + a_3 + \cdots + a_p$, this implies that $a_1 \le N(G) - a_1$. Thus $a_1 \le \frac{N(G)}2$.
\end{proof}

\begin{lemma}\label{thm:Maximum} For any game $G$, we have 
\begin{equation}\label{e:ineq}
N_W(G) \ge \lfloor \log_2 N(G) \rfloor.
\end{equation}
\end{lemma} 

\begin{proof}
We prove this by induction on $N(G)$. As our base case, we consider the position where $N(G) = 1$, when $N_W(G) = 1 > \log_2 N(G) = 0$. Now, we perform the inductive step. Fix a game $G$ and suppose that the result holds for any game $H$ such that $N(H)<N(G)$. Let $n = \lfloor \log_2 N(G) \rfloor$. We divide our analysis into two cases:
\begin{enumerate}
\item If $G$ is a $\mcP$ position, consider a ply $(G, G')$. Then $N(G') \ge \frac{1}{2} N(G)\ge 2^{n-1}$ by Lemma~\ref{l:p2}. Suppose first that Windsor only removes a single candy when going from $G'$ to $G''$, i.e.\ $N(G'') = N(G') - 1$. If $N(G') > 2^{n-1}$, then,
$$N_W(G'') \ge \lfloor \log_2 N(G'') \rfloor \ge n-1,$$
so $$N_W(G) \ge 1 + N_W(G'') \ge n = \lfloor \log_2 N(G) \rfloor,$$ proving the desired result. On the other hand, if $N(G') = 2^{n-1}$, then $N(G'') = 2^{n-1} - 1$ has an odd number of candies and is thus an $\mcN$ position, meaning that Windsor's last ply was invalid.

The only case left to consider is if Windsor removes at least two candies, i.e.\ if $N(G'') \le N(G') - 2$. Since $N(G')-N(G'')>1$, we have
\[N(G')-N(G'')> \lfloor \log_2(N(G')) \rfloor - \lfloor \log_2 (N(G'')) \rfloor.\] 
If $N(G'') = 0$, then $G = [a, a]$ where~(\ref{e:ineq}) holds, and if $N(G'') \neq 0$, then 
\begin{align*} N_W(G)-N_W(G'') &= N(G')-N(G'') \\ &> \lfloor \log_2(N(G')) \rfloor -\lfloor \log_2 (N(G'')) \rfloor \\ &= n-1 - \lfloor \log_2 (N(G'')) \rfloor,\end{align*}
which implies that 
$$N_W(G)> n-1 -\lfloor \log_2 (N(G'')) \rfloor + N_W(G'').$$
Since $N_W(G'') \ge \lfloor \log_2 (N(G'')) \rfloor$, this gives $N_W(G) \ge n$ as desired.
\item Now suppose $G$ is an $\mcN$ position.  Consider a ply $(G, G')$. Since Windsor moves $G \mapsto G'$, $$N_W(G) - N_W(G') = N(G)-N(G')\ge \lfloor \log_2(N(G)) \rfloor - \lfloor \log_2 (N(G')) \rfloor$$ whenever $N(G') > 0$.
By the inductive hypothesis, $N_W(G') \ge \lfloor  \log_2 (N(G')) \rfloor$, so $N_W(G) \ge \lfloor  \log_2 (N(G)) \rfloor$, as desired.
\end{enumerate}
\end{proof}

\begin{lemma}{\label{lem:fairRemove} If $K=[a_1,a_1,a_2,a_2,\ldots,a_p,a_p]$ then for all games $G$, $V(G)=V(G+K)$.}
\end{lemma}
\begin{proof} We prove this by induction on $N(G)+N(K)$. First, the base case $G=K=\varnothing$ is trivial. Now we consider the inductive step. Consider a turn $T = (H := G + K, H', H'')$. If the optimal move is in $G$, then $H' = G' + K$, with $N(G') < N(G)$. Thus, by the inductive hypothesis, $V(H') = V(G' + K) = V(G')$, so
$$V(G+K)=N(G)-N(G')+V(G'+K)=N(G)-N(G')+V(G')=V(G).$$ 
If $(G,G')$ is the optimal ply in $G$, by the same argument we have $$V(G)=N(G)-N(G')+V(G').$$ On the other hand, for any ply $(K,K')$, the opponent can mimic in $K$, and hence move to $H'' = G+K''$ where $K''$ consists of only equal piles. It follows that $$V(G+K')\le N(K)-N(K')+V(G+K'').$$ Thus no move in $K$ can be strictly better than the optimal move in $G$, so we have $V(G+K)=V(G)$, completing the inductive step.
\end{proof}

\begin{lemma} \label{lem:oneremove}For all positive integers $a$, there exists a positive integer $k$ such that $a \oplus (a-1)=2^k-1$.
\end{lemma}
\begin{proof}
If $a$ is odd, then $a \oplus (a-1)$ is $1$, or $2^1-1$. If $a$ is even, then we write $$a = 2^{k_1}+2^{k_2}+\cdots+2^{k_\ell}, \quad k_1>k_2>\cdots>k_\ell>0.$$ Then we have $$a-1 = 2^{k_1}+2^{k_2}+\cdots+2^{k_{\ell-1}}+2^{k_\ell-1}+2^{k_\ell-2}+\cdots+2^{3}+2^2+2^1+1.$$ This gives
\[a\oplus(a-1)=2^{k_\ell}+2^{k_\ell-1}+\cdots+2^3+2^2+2^1+1=2^{k_l+1}-1,\] as desired.
\end{proof}

\begin{lemma} \label{thm:bestarr} The game $G=[1,2,4,8,16,\ldots,2^{n-2},2^{n-1}-1]$ maximizes $N_L(G)$ subject to the constraint that $N(G)=2^{n}-2$. In this case, we have $N_W(G)=n-1$. \end{lemma}

\begin{proof} 
First, let us compute $N_W([1,2,4,8,16, \ldots ,2^{n-2},2^{n-1}-1])$. If Luca removes the entire largest pile, then Windsor is forced to remove a single candy, leaving the game $G'=[1,2,4,8,16,\ldots,2^{n-3},2^{n-2}-1]$. When $n=2$ we get $N_W=1$. By induction $N_W = n-1$. By Lemma~\ref{thm:Maximum} $N_W(G) \ge n-1$, for an arbitrary $\mcP$ position $G$ with $N(G)=2^n-2$. Since $G=[1,2,4,8,\ldots,2^{n-2},2^{n-1}-1]$ achieves equality, it minimizes $N_L$ subject to $N(G)=2^n-2$.
\end{proof}

\begin{lemma} \label{thm:Strat}
Given \[G = [1,2,4,8,\ldots,2^{k-2},2^{k},\ldots,2^{n-2},2^{n-1}-1-2^{k-1}],\] the ply $P = (G, G')$ with \[G' = [1,2,4,8,\ldots,2^{k-2},2^{k},\ldots,2^{n-2},2^{k-1}]\] is an optimal move. Then, \[N_W(G) = N_W(G') = n -1.\] 
\end{lemma} 
\begin{proof} Lemma~\ref{thm:Maximum} shows that it is impossible for Luca to concede fewer than $n-1$ candies to Windsor. Therefore to show optimality, it suffices to show that $N_W(G) = n - 1$. If Luca moves the pile of size $2^{n-1}-1-2^{k-1}$ to a pile of size $2^{k-1}$, the remaining game $G'$ has $2^n - 2$ candies, with $N_W(G') = n - 1$ by Lemma~\ref{thm:bestarr} assuming optimal play by Luca. Thus, $G$ minimizes $N_W$ with $N_W(G) = n-1$ as desired.
\end{proof}

\begin{thm} \label{thm:equalitycases}
Given a game $G$, $N_W(G) \ge \lfloor \log_2(N(G)) \rfloor$. Equality is achieved only when $N(G) = 2^n, 2^n - 2$, or $2^n - 2^k - 2$, for $n, k \in \ZZ^+$, $n > k + 1,n>2$ in the following arrangements:
\begin{enumerate}
\item $N(G) = 2^n$ and $G = [1,1,1,2,4,8,\ldots,2^{n-2},2^{n-1}-1]$
\item $N(G) = 2^n - 2$ and $G = [1,2,4,8,\ldots,2^{n-2},2^{n-1}-1]$
\item $N(G) = 2^n - 2^k - 2$ and $G = [1,2,4,8,\ldots,2^{k-2},2^k, \ldots ,2^{n-2},2^{n-1}-1-2^{k-1}]$
\end{enumerate}
\end{thm}

\begin{proof} First, note that it is sufficient to prove the result when $G$ is a $\mcP$ position. To see this, suppose that we have proven the result for all $\mcP$ positions, and $G$ is an $\mcN$ position with a ply $(G,G')$ where $G'$ is a $\mcP$ position, then we have \[N_W(G)\ge N(G)-N(G')+N_W(G')\ge N(G)-N(G')+\lfloor\log_2(N(G'))\rfloor\ge\lfloor\log_2(N(G))\rfloor.\] Thus from now on, we shall always assume that $G$ is a $\mcP$ position.

We prove the result by induction on $N(G)$. Via a finite check, this is true whenever $N(G) \le 16$. For the inductive step, suppose that equality is achieved only in the above positions for all positions with $N(G) < M$. We want to show that if $N(G)=M$, this theorem holds. 

First, we show that $N_W(G)=\lfloor \log_2 (N(G)) \rfloor$ implies $N_W(G'')=\lfloor \log_2(N(G'')) \rfloor$. Let $M = 2^n + x$ where $n = \lfloor \log_2 M \rfloor$. Then, $$2^{n-1} \le 2^{n-1} + \frac x2 \le N(G') \le 2^n + x -1 \le 2^{n+1}.$$
If $N(G')- N(G'')=1$, then $$2^{n-1}-1 \le N(G'') < 2^{n+1}-1.$$ Since $G$ is a $\mcP$ position, $N(G)$ is even and thus $N(G'') \neq 2^{n-1}-1$. Thus $2^{n-1} \le N(G'') < 2^{n+1}-1$, so by Lemma~\ref{thm:Maximum}, $N_W(G'') \ge n-1$. If
$N_W(G'') \ge n$, then $N_W(G) \ge n+1$, so in any potential equality case, we must have $N_W(G'') = n-1$. If $N(G')- N(G'') \ge 2$, then whenever $ N(G'') > 0$ and $N(G')-N(G'')>1,$ we have \[ N(G')-N(G'')> \lfloor \log_2(N(G')) \rfloor - \lfloor \log_2 (N(G'')) \rfloor.\]
Since $N(G'') \ge \lfloor \log_2 (N(G'')) \rfloor$,  if $N_W(G) = \lfloor \log_2 (N(G)) \rfloor$, then $N_W(G'') = \lfloor \log_2 (N(G'')) \rfloor$. If $N(G'')=0$, then $N(G') = \frac{N(G)}2$.
Thus, $N_W(G)=\lfloor \log_2 (N(G)) \rfloor$ which implies $N_W(G'')=\lfloor \log_2(N(G'')) \rfloor$. Therefore, by the inductive hypothesis, $G''$ must one of the three positions above.

Now we show that if $G''$ is any one of the above three positions, then so is $G$, thereby completing the induction.
\begin{enumerate}
\item 
Suppose that $$G''=[1,1,1,2,4,8, \ldots, 2^{n-3},2^{n-2}-1].$$ In order to have $$N_W(G')-N_W(G'') = \lfloor \log_2 (N(G)) \rfloor - \lfloor \log_2 (N(G'')) \rfloor,$$ we must have $N(G) \in\{2^n, 2^n+2\}$. If $N(G)=2^n+2$, then $N(G)-N(G')=2^{n-1}+1$. However, this implies that $$G = [2^{n-1}+1,2^{n-1}+1]\text{ or }[1,2^{n-1},2^{n-1}+1],$$ since those are the only two $\mcP$ positions with $N(G)=2^n+2$. Neither of those can produce a $G''$ of the specified form. Therefore, $N(G)=2^n$ and $N(G)-N(G')=2^{n-1}-1$. and there must have been a pile of size at least $2^{n-1}-1$ in $G$. If there was a pile of size at least $2^{n-1}$, we have the same issue as above with $2^n + 2$. Consequently, there must be a pile of size exactly $2^{n-1}-1$. If $N(G)=2^n$, Windsor removed $1$ candy on the first term, giving the Grundy value of $G'$, $\mcG(G') \in \{1,3, 2^{n-1}-1\}$. In the first two cases, there is no way to achieve $N(G)-N(G')=2^{n-1}-1$. Therefore, 
$$G=[1,1,1,2,4,8,\ldots,2^{n-2},2^{n-1}-1].$$
\item Now suppose that $$G''=[1,2,4,8,\ldots,2^{n-3},2^{n-2}-1].$$ As Windsor removed $1$ candy, $\mcG(G') \in \{1, 3, 2^{n-1} - 1\}$. If $\mcG(G') = 1$, then
$$G = [1,2,4,8,\ldots,2^{l}+1,\ldots,2^{m}+1,\ldots,2^{n-3},2^{n-2}-1],$$ which allows Windsor to remove $1$ candy from a different pile to increase his winnings, contradicting optimal play. If $\mcG(G') = 3$, then $$G = [2,2,4,8,\ldots,2^{l}+3,\ldots,2^{n-2}-1]\text{ or }[2,2,3,4,\ldots,2^{n-2}-1].$$ Windsor could have removed $3$ from the $2^{n-2}-1$ and received more candies while still winning, again contradicting optimal play. If $\mcG(G') = 2^{n-1}-1$, then $$G'=[1,2,4,8, \ldots, 2^{n-3},2^{n-2}].$$ So, either Luca moved from $2^{n-1}-2^k-1$ to $2^k$ or from $2^{n-1}-1$ to $0$. The first case gives the third game above, and the second gives the second game above.
\item Finally, suppose that $$G''=[1,2,4,8,\ldots,2^{k-2},2^{k},\ldots,2^{n-3},2^{n-2}-2^{k-1}-1]$$ with $\mcG(G') \in \{1, 3, 2^{k}-1\}$. If $k \ge 2$, then as $G-G''  \ge 2^k+1$, it would be impossible for Windsor to both remove one, and have $\lfloor\log_2(G'')\rfloor < \lfloor \log_2(G) \rfloor $. Otherwise $k = 1$, $\mcG(G') = 1$, and thus $G = 2 + G''$ so $\lfloor \log_2(G'') \rfloor = \lfloor \log_2 (G) \rfloor$, a contradiction.

\end{enumerate}
 \end{proof}

\begin{thm} \label{thm:5pile} For all $N \in \ZZ^+$, there exists a 5-pile game $G$ with  $N(G) = N$ and $N_W(G) \leq \frac{3}{2}\sqrt{2N} -2 $.
\end{thm} 

\begin{proof} 
We can write $N$ in binary as
\[ N = 2^{k_1}+ 2^{k_2}+\cdots+2^{k_n}+2^{k_{n+1}}+2^{k_{n+2}}+ \cdots+2^{k_{n+m}},\] where $k_1>k_2>\cdots>k_{n+p}$, where $n$ is defined so that $k_n\ge\lfloor\frac{k_1}{2}\rfloor$ but $k_{n+1} < \lfloor\frac{k_1}{2}\rfloor$. Thus $n$ is the minimal $i$ such that $2^{k_{i+1}} < \sqrt{N}$. 
Consider the game $G_1 = \mfG(a,m,0)$, where
\[m=2^{k_1-\lfloor \frac{k_1}2 \rfloor}+2^{k_2-\lfloor \frac{k_1}2 \rfloor}+2^{k_3-\lfloor \frac{k_1}2 \rfloor}+\cdots+2^{k_n-\lfloor \frac{k_1}2 \rfloor}-1\qquad\text{and}\qquad
a=2^{\lfloor \frac{k_1}2 \rfloor-1}-1.\]
By construction, $N(G_1) < N$. From this, we can construct the game
\[ G = \left [2^{\lfloor \frac{k_1}2 \rfloor-1}-1,2^{k_1-1}+2^{k_2-1}+ \cdots +2^{k_n-1}-2^{\lfloor \frac{k_1}2 \rfloor-1}, \right .2^{k_1-1}+2^{k_2-1}+ \cdots \]\[ \left . \cdots +2^{k_n-1}-1,2^{k_{n+1}-1}+2^{k_{n+2}-1}+\cdots+2^{k_{n+m}-1} -1,2^{k_{n+1}-1}+2^{k_{n+2}-1} + \cdots +2^{k_{n+m}-1}-1 \right ] \] 
where $N(G) = N$. Note that the last two piles of $G$ are identical.
Corollary~\ref{cor:UsefulBound} gives
\[N_W(G_1) \le 2^{\lfloor \frac{k_1}2 \rfloor}-2+2^{k_1}+2^{k_2}+ \cdots +2^{k_n}-2^{\lfloor \frac{k_1}2 \rfloor}-
(2^{k_1}+2^{k_2}+ \cdots +2^{k_n}) + r_N, \,\, r_N \le \sqrt{2N}, \]
and therefore \[N_W(G) \le \frac32 \sqrt{2N} -2\]
because \[N_W(G)=2^{k_{n+1}-1}+2^{k_{n+2}-1}+\cdots+2^{k_{n+p}-1}+r_N-2.\] 
\end{proof}

\begin{thm} \label{thm:k-1} If $G$ is a game containing $p$ piles with no duplicate piles, then $N_W(G) \ge p-1$.
\end{thm}

\begin{proof}
We shall prove this by induction on $N(G)$. When $N(G)<2$, the result is trivial. When $N(G)=2$, then $G$ is either $[2]$ or $[1,1]$. The first one gives $2 \ge 0$, and the second $1 \ge 1$, as desired. Suppose the claim is true for all $G$ with $N(G)<n$. We show it holds when $N(G) = n$.

\begin{itemize} \item[$(\mcN)$] Let $G$ be an $\mcN$ position. Then the number of piles in $G'$ is at most one fewer than that of $G$. So, either $N_W(G') \ge p-2$ or Windsor made a move to make two piles equal sizes. In the first case, Windsor must have removed at least one candy, so $N_W(G)\ge p-1$ as desired. If Windsor moved to create a duplicate pile, $$G' = [a,a,g_1,g_2,g_3,\ldots,g_{p-2}],$$  where the $g_i$'s are all distinct. By \ref{lem:fairRemove}, $N_W(G')=a + N_W([g_1,g_2,\ldots,g_{p-2}])$. By induction $N_W([g_1,g_2,\ldots,g_{p-2}]) \ge p-3$. As $a \ge 1$, we get that $N_W(G') \ge p-2$ so $N_W(G)\ge p-1$ as desired 

\item[$(\mcP)$] Suppose $G$ is a $\mcP$ position. 
\begin{enumerate}
\item If Luca doesn't remove a fill pile, then $G'$ has the same number of piles as $G$. 
We consider cases:
\begin{enumerate}
\item If there are no duplicates in $G'$, by the inductive hypothesis, $N_W(G)=N_W(G') \ge p-1$ as desired.
\item Suppose Luca creates a duplicate pile, so $$G' = [a,a,g_3,\ldots,g_p].$$ Then we have $N_W(G')=a+N_W([g_3,\ldots,g_p])$. If $a \neq 1$, $N_W(G)=N_W(G') \ge 2+ p-3 =p-1$, via inductive hypothesis. Suppose $$G'=[1,1,g_3,\ldots,g_p].$$ In that case, we must have had $G=[g_1,g_2,g_3,\ldots,g_p]$ with $g_1 = 1$. Windsor cannot move in a $1$ pile, or else Luca would have been able to move to $G'' = [g_2,g_3, \ldots, g_p]$, contradicting the assumption that $G\in\mcP$. So, his winning move must be in one of the piles $g_3,\ldots,g_p$. If Windsor doesn't remove a pile, we get 
\begin{align} \label{eq:NWPpos}  \begin{split}N_W(G) &=N_W(G') \\ &= 1+N_W([g_3,g_4,\ldots,g_p]) \\ &= 2+N_W([z,g_4,\ldots,g_p]) \\ &\ge 2 + p - 3 \\ &= p-1.\end{split} \end{align} The first equality in~(\ref{eq:NWPpos}) is because it is currently Luca's move. The second equality follows from Lemma~\ref{lem:fairRemove}, The third equality follows because Windsor removed one candy. The fourth inequality follows from the inductive hypothesis. 
\item If Luca first creates a $1,1$ duplicate (i.e. moves a pile $g_2$ to size $1$ with an existing pile $g_1$ of size $1$) to obtain $G'$, Windsor removes a pile in $G'$. We have \begin{align*} N_W(G) &= N_W(G') \\ &= 1+N_W([g_3,g_4,\ldots,g_p] \\ &= 1+g_3+N_W([g_4,\ldots,g_p]) \\ &\ge 1+g_3+p-4,\end{align*} where $g_3$ is the pile Windsor removes. If $g_3 \neq 1$, we have $$N_W(G) \ge 1+2 +p-4 =p-1,$$ as desired. But if $g_3=1$, then $G$ had a $1,1$ duplicate already, contrary to hypothesis.
\end{enumerate}

\item Suppose Luca removes a pile. We have $G' = [0,g_2,g_3,\ldots,g_p]$.
We further subdivide into cases:
\begin{enumerate}
\item If Windsor removes a pile $g_2$, then it is Luca's turn so $g_1 \oplus g_2 = 0$ and $g_1 = g_2$, giving an initial duplicate pile.
\item Suppose Windsor doesn't remove a pile and creates no duplicate piles when he moves $G'$ to $G''$. Via the inductive hypothesis $N_W(G'') \ge p-2$. Since Windsor removed at least one candy, $N_W(G) \ge p-1$ as desired.
\item Suppose Windsor removes no entire pile, but creates some duplicate pile of size $a \ge 2$ so $G'' = [a, a, g_4, g_5, \ldots g_p]$ with $$N_W(G'')=a+N_W([g_4,g_5,\ldots,g_p])\ge a+p-4.$$
Since $a \ge 2$, and Windsor removed at least one candy, $$N_W(G)\ge 1+ a+ p -4 \ge 1+2 + p-3=p-1$$ as desired. 
\item Finally, suppose that Windsor removes some candies to create a duplicate pile of size $1$ with $G'' =[1, 1, g_4, g_5 , \ldots, g_p]$. This would give
$$G' = [1,2,g_4,g_5,\ldots,g_p], \quad G = [1,2,3,g_4,g_5,\ldots,g_p]$$
as Luca removed a pile (so no other pile had size $2$).
Since $G'' \in \mcP$, $\mcG(G') = 3$. It suffices to show that if $H = [g_4, g_5, \ldots, g_p]$, that $N_W(H) \ge p-3$. If $H = \varnothing$, we are done. Thus suppose $H$ has at least one pile.
Note that for all $i \ge 4$, $g_i > 1$ and $g_i \equiv 0, 1 \pmod{4}$, and all these piles of $H$ are distinct. We can consider the possible moves in Luca's ply $(H, H')$ as we did above. 
\begin{itemize}
\item[$\bullet$] Any duplicate pile created has size at least $4$, so creating a duplicate pile would yield the desired bound: \[N_W(H)=N_W(H') \ge 4+N_W([g_6,g_7,\ldots,g_p]) \ge 4+p-6 =p-2.\]
\item[$\bullet$] If Luca neither removes a pile nor creates a duplicate, Windsor must move in a distinct pile from Luca. If Windsor removed a pile, he removed at least $4$ candies, so $N_W(G) \ge n-1$. Since $H$ is duplicate free, Windsor cannot create a duplicate. Thus, if Windsor didn't remove a pile, we thus obtain $H'' = [a, b, g_6, g_7, \ldots, g_p]$ with
\begin{align*} N_W(G) &= N_W(G') \\ &\ge 1+N_W(G'') \\ &= 2+N_W(H) \\ &= 2+N_W(H') \\ &= 3+N_W([a,b,g_6,g_7,\ldots,g_p]) \\ &\ge 3+n-4 \\ &= n-1,\end{align*}
\item[$\bullet$] Suppose Luca removes a pile. Since $g_i \equiv 0, 1 \pmod 4$ for all piles in $H$, Windsor must have removed at least $3$ candies since $H'' \in \mcP$; furthermore, because $H$ contains no duplicates, Windsor cannot have removed an entire pile in moving from $H'$ to $H''$. Thus $H''$ consists of $p-4$ piles. If $H''$ has no duplicate piles, then by induction $N_W(H'')\ge p-5$, so $$N_W(H)\ge 3+(p-5)=p-2,$$ which is greater than the required $p-3$. On the other hand, if $H''$ has a duplicate pile, say with $H''=[g_6,g_6,g_7,g_8,\ldots,g_p]$, then \begin{align*} N_W(H) &\ge 3+N_W(H'') \\ &= 3+g_6+N_W([g_7,g_8,\ldots,g_p]) \\ &\ge 3+g_6+(p-6) \\ &\ge p-3.\end{align*}
\end{itemize}
 This completes the analysis of the final case, and shows the desired inductive hypothesis. 
\end{enumerate}
\end{enumerate}
\end{itemize}

\end{proof}

For small $N$, we can use the above results to identify the games $G$ with $N(G) = N$ that minimize $N_W(G)$.

\begin{ex} \label{lem:bestten}
For $N = 10, 12, 16$ we compute the unique games $G$ that minimize $N_W(G)$ via Theorem~\ref{thm:equalitycases}.
\begin{itemize}
\item If $N=10$, then $G = [1, 4, 5]$ minimizes $N_W$, with $N_W(G) = 3$.
\item If $N=12$, then $G=[2,4,6]$ minimizes $N_W$, with $N_W(G) = 3$.
\item If $N = 14$, then $G=[1,2,4,7]$ minimizes $N_W$, with $N_W(G) = 3$.
\item If $N=16$, then $G=[1,1,1,2,4,7]$ minimizes $N_W$, with $N_W(G) = 4$.
\end{itemize}
\end{ex}

\section{Conjectures and Concluding Remarks} \label{sec:examples}

\subsection{$4$-Pile \textsc{Candy Nim}}
Most of our attention with respect to strategies and bounds on $V(G)$ has been focused on the case when $G$ is a $3$ pile game. We include a brief analysis and several conjectures regarding $V(G)$ and optimal play for $4$ pile games.

First, we show that, in the 4-pile game, Luca does not always have an optimal move in the largest pile.

\begin{ex} Let $G=[1,5,16,20]$. We show $V([1,5,16,20])=28$, where Luca's optimal move is to remove three candies from the pile of size 5. By checking, we have the following optimal game play:
$$[1, \color{red}5\color{black}, 16, 20] \overset{\textcolor{red}{L}}\rightarrow [1, \color{red}2\color{black}, 16, \color{green}20\color{black}] \overset{\textcolor{green}{W}}\rightarrow [1, 2, 16, \color{green}1 \color{red}9\color{black}] \overset{\textcolor{red}{L}}\rightarrow [1,2,\color{red}12, \color{green}16\color{black}] \overset{\textcolor{green}{W}}\rightarrow [1,2,12,\color{green}1\color{red}5\color{black}] \overset{\textcolor{red}{L}}\rightarrow [1,2,\color{red}8,\color{green}12\color{black}]$$ 
$$\overset{\textcolor{green}{W}}\rightarrow [1,2,8,\color{green}1\color{red}1\color{black}]\overset{\textcolor{red}{L}}\rightarrow [1,2,\color{green}8,\color{red}4\color{black}] \overset{\textcolor{green}{W}}\rightarrow [1,2,\color{green}7\color{black},4] = [1, 2, 4, 7]$$
By Theorem~\ref{thm:equalitycases}, $V([1, 2, 4, 7]) = 8$. Thus,
$$V(G) = 20 + 8 = 28$$
\end{ex}

We can obtain lower bounds on some families of four pile games $G$ using related three pile games. We first consider the four pile games $G$ with smallest two piles of size $1, 2$, and show that their values $V(G)$ are bounded by the ``corresponding'' $3$-pile game with smallest pile size $3$.

\begin{prop} \label{p:4pile}
Let $m$ be a positive integer. Then both of the following hold:
$$V([1,2,4m, 4m + 3]) \geq V([3,4m, 4m+3]),$$ $$V([1,2,4m+1, 4m+2]) \geq V([3,4m+1,4m+2]).$$
\end{prop}

\begin{proof}
Let \begin{align*} G_1 = [3,4m, 4m+3], && G_2 = [3,4m+1, 4m+2], \\ H_1 = [1,2,4m,4m+3], && H_2 = [1,2,4m+1, 4m+2].\end{align*}
We will show the desired result by induction on $m$. Let $m = 1$ be our base case.
We can check $6 = V([3,4,7]) \leq V([1,2,4,7]) = 8$ and $V([3,5,6]) = V([1,2,5,6]) = 6$.
Given $i\in\{1,2\}$, for every possible optimal turn $T_{G_i} = (G_i,G_i',G_i'')$ we must show that there exists a turn $T_{H_i} = (H_i, H_i', H_i'')$ such that $$V_{T_{G_i}}(G_i) + V(G_i'') \leq V_{T_{H_i}}(H_i) + V(H_i'').$$

Suppose that $V(G_i) \leq V(H_i)$ for $m<w$.
Let $m = w$. Suppose $G_i' = [3,a,b]$ and $G_i'' = [3,a,c]$. Then we set $H_i' = [1,2,a,b]$ and $H_i'' = [1,2,a,c]$, so that $V_{T_{G_i}}(G_i) = V_{T_{H_i}}(H_i)$ and $V(G_i'') \leq V(H_i'')$ by the inductive hypothesis. If $G_i' = [2,a,b]$ or $G_i' = [1,a,b]$, then we set $H_i' = [0,2,a,b]$ and $H_i' = [1,0,a,b]$, respectively. This yields $$V_{T_{G_i}}(G_i) + V(G_i'') = V_{T_{H_i}}(H_i) + V(H_i'').$$

Now suppose that $G_i'' = [0,a,b]$. If $i = 1$, then $V_{T_{G_1}}(G_1) + V(G_1'') = 0$ and $V_{T_{H_1}}(H_1) + V(H_1'') \geq 0$ by Lemma \ref{lem:valuegeq0}.
If $i=2$, then $V_{T_{G_2}}(G_2) + V(G_2'') \leq 2$, which implies that it is not an optimal move since Luca could instead remove the largest pile in $G_2$ and obtain an overall value of $4$. 
Thus, by induction, we have $$V([1,2,4m, 4m + 3]) \geq V([3,4m, 4m+3]), \quad V([1,2,4m+1, 4m+2]) \geq V([3,4m+1,4m+2]).$$
\end{proof}

\subsection{General Play}
We can hope to make even more general inferences from the $3$-pile game to multi-pile \textsc{Candy Nim} games. Notably, we conjecture a similar result to Proposition~\ref{p:4pile} holds for a broader family of \textsc{Candy Nim} games.

\begin{conj}\label{c:gen}
Suppose $G = [a, b, c]$ with $a < b < c$. Then for some $j > 1$, there exist $a_1, \ldots, a_j$ with 
$$a = a_1 + \cdots + a_j = a_1 \oplus \cdots \oplus a_j$$
such that the game $$H = [a_1, a_2, \ldots, a_j, b, c]$$ satisfies $V(H) \geq V(G)$.
\end{conj}

\begin{rem} 
Note that it is \textit{not true} that for any game $G =[a, b,c]$, all such decompositions of $a = a_1 + \cdots + a_j = a_1 \oplus \cdots a_j$, with resulting game $H$ as in Conjecture~\ref{c:gen} satisfies $V(H) \ge V(G)$. As a counterexample, consider the game $G = [31,42,53]$, with $a = 31$. Using the decomposition $a_1 = 1, a_2 = 2, a_3 = 4, a_4 = 8, a_5 = 16$ we obtain the game $H = [1,2,4,8,16,42,53]$. However, $V(G) = 96$ while $V(H) = 94$.
\end{rem}

We can also hope to extend the analysis of Section~\ref{sec:multipile}. Observationally, for a fixed number of candies, the games $G$ that optimize $N_W(G)$ have specific structural properties that we conjecture hold in general:

\begin{conj}
For all fixed $N > 0$, there exist (not necessarily distinct) games $G_1, G_2$ with $N(G_1) = N(G_2) = N$ such that 
$$N_W(G_1) = N_W(G_2) = \max_{H;\,N(H) = N} N_W(H)$$
where $G_1$ has a pile with at least $N/4$ candies and $G_2$ has at most $c \log N$ piles for some absolute constant $c > 0$.
\end{conj}

\bibliography{candynim}
\bibliographystyle{alpha}

\end{document}